\documentclass[10pt]{elsarticle}
\usepackage{amssymb}
\usepackage{amsmath}
\usepackage{amsfonts} 
\usepackage{amsbsy}
\usepackage{amscd}
\usepackage{amsthm}
\usepackage{courier}
\usepackage{fullpage}
\usepackage{times}
\usepackage{psfrag}
\usepackage{graphics}
\usepackage{color}
\usepackage{bbm}
\usepackage{colortbl}
\usepackage{graphicx}
\usepackage{caption}
\usepackage{bbm}
\usepackage{subcaption}
\usepackage{array,supertabular}
\usepackage{tabularx}
\usepackage{booktabs}
\usepackage{multirow}
\usepackage{ulem}
\usepackage{units}
\usepackage{mathabx}
\usepackage{accents}
 \usepackage{relsize}
\usepackage{xfrac}
\usepackage{algpseudocode}
\usepackage{rotating}
\usepackage[capitalize]{cleveref}
\usepackage{mathrsfs} %better calligraphic fonts (\mathscr instead of \mathcal)
\usepackage[framemethod=tikz]{mdframed} %%to make a framed environment
\usepackage{hypbmsec}
\usepackage[framemethod=tikz]{mdframed}
\usepackage[most]{tcolorbox}
\usepackage{lipsum}% just to generate text for the example

\usepackage{appendix}
\usepackage{float}
\usepackage{makecell}

\crefname{secinapp}{Section}{Sections}
\Crefname{secinapp}{Section}{Sections}

\newlength{\dhatheight}

\newtheorem{theorem}{Theorem}[section]

\newtheorem{example}[theorem]{Example}

\theoremstyle{definition}
 
\theoremstyle{remark}
\newmdtheoremenv[
hidealllines=true,
leftline=true,
innertopmargin=0pt,
innerbottommargin=0pt,
linewidth=4pt,
linecolor=gray!40,
innerrightmargin=0pt,
innertopmargin=-6pt,
]{examplei}{Example}

\begin{document}
%\fancyhead{}

\begin{frontmatter}
\title{Sums of Squares: Methods for Proving Identity Families}
 \author[tu]{Russell Jay Hendel \\ RHendel@Towson.Edu}
\address[tu]{Dept of Mathematics, Towson University, Towson Maryland 21252, USA}

\begin{abstract}
This paper presents both a result and a method. The result presents
a closed formula for the sum of the first $m+1,m \ge 0,$ squares of the 
sequence $F^{(k)}$ where
each member is the sum of the previous $k$ members and with initial conditions
of $k-1$ zeroes followed by a 1. The generalized result includes the known result
of sums of squares of the Fibonacci numbers and recent results of 
Ohtsuka-Jakubczyk, Howard-Cooper, Schumacher, and Prodinger-Selkirk for the cases $k=2,3,4,5,6.$
The paper contributes a closed formula for coefficients for all $k.$
To prove the result, the paper introduces  a new method,
  the \emph{algebraic verification} method, which reduces proof of an identity
to verification of the equality of finitely many pairs of finite-degree polynomials, 
possibly in several variables. 
Additionally, the paper provides a visual aid, \emph{labeled index squares},
for complicated proofs. Several other papers proving families of identities are examined; it is 
suggested that the collection of the uniform proof methods used 
in these papers could possibly produce a new 
trend in stating and proving identities.

\end{abstract}

\begin{keyword}
$k$-bonacci, sums of squares, family of identities, verification, 
generalized Fibonacci,
\end{keyword}
 
\end{frontmatter}

%-------------------------------
% SEE CUT OUT MATERIAL

\section{Motivation}\label{sI}

 The generalized Fibonacci numbers, also known as the $k$-bonacci numbers, are
defined by 
\begin{equation}\label{recursion}
	F^{(k)}_n= \sum_{i=1}^k F_{n-i},\qquad \text{ with } 
	F_i=0, \text{ for $0 \le i \le k-2,  F_{k-1}=1,$} k\ge 2.
\end{equation}
These  particular initial conditions
  are consistent with the Online Encyclopedia of 
Integer Sequences (OEIS). For $k=2,3,4,5$ we obtain the Fibonacci \cite[A000045]{OEIS},
Tribonacci \cite[A000073]{OEIS}, Tetranacci  \cite[A000078]{OEIS}, and
Pentenacci numbers \cite[A001591]{OEIS} respectively.

Closed formulas for the sum of the first $m+1$ squares of $\{F^{(k)}_n\}_{n \ge  0},$
\begin{equation}\label{TheProblem}
	\displaystyle \sum_{i=0}^m F^{(k) \; 2}_i,
\end{equation}
are known for the cases $k=2,$ \cite[pp. 77--78]{Koshy}, 
$k=3,$ \cite{Jakubczyk,Ohtsuka,Schumacher3},
$k=4,$ \cite{Jakubczyk, Prodinger_Selkirk, Schumacher4}, $k=5,$ \cite{Prodinger_Selkirk},
and
$k=6,$ \cite{Jakubczyk}.

For $k=2,$ \cite[pp. 77--78]{Koshy}, the closed formula for \eqref{TheProblem} is
$\sum_{i=0}^m F_i^2 = F_m F_{m+1}.$ 
 \cite{Howard_Cooper} shows that for all $k,$ the closed formula for \eqref{TheProblem} has a
$F^{(k)}_m F^{(k)}_{m+1}$ term, generalizing the closed formula for the case $k=2.$ The algebraic
methods of \cite{Schumacher4} and the generating-function methods of
\cite{Prodinger_Selkirk}, would appear to be able to provide closed formulas for 
\eqref{TheProblem} for arbitrary $k;$ however, this is not done explicitly in these papers.
Moreover, neither of these 
methods explicitly identify the patterns in the numerical coefficients of the closed formulas 
for general
$k.$

The main contributions of this paper are i) an explicit formula for the numerical coefficients of 
the closed formula for general $k,$ ii) the \emph{algebraic verification} method which 
provides a uniform proof of the closed formula for all $k,$ and iii) the 
\emph{labeled index square}
method which provides visual aids facilitating proofs.  

Since the algebraic verification method is new, we should compare it to other methods.
First, the algebraic verification method can simply be regarded as another tool
to prove identities along with the Binet form, generating functions, and matrices.
It differs from the generating-function method 
\cite{Prodinger_Selkirk, Schumacher3} in that the proof has no algebraic prerequisites;
for example, 
\cite{Prodinger_Selkirk} succeeded in their paper by using an established 
theory for the Binet form of generalized Fibonacci numbers.

Both the generating functions method and algebraic verification are typically computational.
A simple glance at \cite{Schumacher3} shows many lines of algebraic manipulations. 
In \cite{Prodinger_Selkirk} these manipulations are done by software. The algebraic verification
method used in this paper requires computation of the coefficients of 24 polynomials. However, while
there is a lot of work, each of the underlying polynomials
containing these coefficients are polynomials in  at most 3 variables
and are polynomials of at most degree 2; 
the computations are done by pencil and paper without a need for software.

In comparing the algebraic verification method with other algebraic-recursive methods,
we see that the algebraic verification method is more direct. Once some routine simplifications
are done, it is immediately transparent what has to be checked to prove the theorem. Contrastively,
the methods of \cite{Howard_Cooper, Jakubczyk,Ohtsuka,Schumacher4} all require some clever algebraic manipulation to accomplish the proof.

Another important difference between the methods is that for example \cite{Prodinger_Selkirk}
gives methods of computing coefficients for closed formulas for sums of squares of generalized
Fibonacci numbers, without showing that these coefficients have a pattern expressible as a 
polynomial of at most degree 2 in three variables. The identification of these patterns, as 
mentioned above, is one contribution of this paper.
  
The idea of proof methods that \emph{uniformly} prove families of identities on recursive
families has independent interest beyond the proof of the particular result in this paper. 
It also suggests a new trend for stating and proving identities.

These ideas motivate the following outline to this paper. 
Section \ref{sE} presents the closed
formula for sums of the first $m+1$ squares of the $F^{(k)}$ for $k=2,3,4,5,6.$
 Tables of these coefficients exhibit patterns of regularity motivating the statement
 of the Main Theorem in Section \ref{sMT} which also explains the core idea of the algebraic
verification method. Section \ref{s7} introduces labeled index-squares, a visual aid
facilitating complex proofs. The following section presents a complete outline of the proof with an illustrative example. 
This is followed by two sections, 
which break up the sum of squares in the 
statement of the Main Theorem into five groups and analyze their coefficient patterns
using the labeled index squares.
Section \ref{sDistribution} then 
combines the results of previous
sections and completes the proof. Section \ref{sC} concludes the paper with a review of
several recent results on families of identities and speculates on a possible new trend
in approaching  identities.
 
\section{Examples}\label{sE}

This section presents closed formulas for \eqref{TheProblem} for $k=2,3,4,5,6.$
Tables created from these formula exhibit patterns which will enable us to formulate the Main Theorem
in the next section.

First we present the examples. References for closed formulas
for the cases $k=2,3,4,5,6$ were presented in the introductory section.
 Note, that  in these examples,
and for the rest of the paper, we ease notation by letting $G$ stand for $F^{(k)}.$

\begin{equation}\label{ex2}
\mathbf{For \;  G_n=F^{(2)}_n:  }	\sum_{i=0}^m G_i^2 = 
			\frac{1}{2} \biggl( 2G_m G_{m+1} \biggr) =G_m G_{m+1}. 
\end{equation}

\begin{equation}\label{ex3} 
\mathbf{For \; G_n=F^{(3)}_n:}	\sum_{i=0}^m G_i^2 = 
	\frac{1}{4}\biggl(-G_m^2 -4 G_{m+1}^2-G_{m+2}^2 +
		2G_m G_{m+2} + 4G_{m+1} G_{m+2} +1\biggr). 
\end{equation}

\begin{multline}\label{ex4}
 \mathbf{For \; G_n=F^{(4)}_n:}
 \displaystyle \sum_{i=0}^{m} G_i^2 =\frac{1}{6} 
	\biggl(-2G_m^2 -8 G_{m+1}^2-6 G_{m+2}^2  
	- 2G_{m+3}^2
	-2 G_m G_{m+1} \\
	+ 2G_m G_{m+3} + 4 G_{m+1} G_{m+3} +6 G_{m+2} G_{m+3} + 2\biggr).
\end{multline}
 
\begin{multline}\label{ex5}
\mathbf{For \; G_n=F^{(5)}_n:}
  \displaystyle \sum_{i=0}^{m} G_i^2 = \frac{1}{8} 
	\biggl(-3 G_m^2 -12 G_{m+1}^2- 11 G_{m+2}^2 - 8 G_{m+3}^2 - 3 G_{m+4}^2  \\
		-4 G_m G_{m+1} -2 G_m G_{m+2}  +2 G_m G_{m+4} 
	-4 G_{m+1} G_{m+2}\\ +4 G_{m+1} G_{m+4} 
	 +6 G_{m+2} G_{m+4} +8 G_{m+3} G_{m+4} +3\biggr).
\end{multline}
 
\begin{multline}\label{ex6}
\mathbf{For \; G_n=F^{(6)}_n:}
  \displaystyle \sum_{i=0}^{m} G_i^2 = \frac{1}{10} 
	\biggl(-4 G_m^2 -16 G_{m+1}^2- 16 G_{m+2}^2 - 14 G_{m+3}^2 - 10 G_{m+4}^2 -4 G_{m+5}^2 \\ 
		-6 G_m G_{m+1} -4 G_m G_{m+2}  -2 G_m G_{m+3} +2 G_{m} G_{m+5}
		-8 G_{m+1} G_{m+2}   -4 G_{m+1} G_{m+3}\\  +4 G_{m+1} G_{m+5}     
		-6G_{m+2} G_{m+3}  +6 G_{m+2} G_{m+5}  +8 G_{m+3} G_{m+5} 
		+10 G_{m+4} G_{m+5} + 4 \biggr).
\end{multline}

We next discuss identifying the patterns in the coefficients in these identities and the consequent closed formulas. The general summand on the right-hand side of the above equations is
$$
	\frac{N_{i,j}}{D_k} G_{m+i} G_{m+j};
$$
there is also a constant term, the last summand on the right-hand side.
(To ease notation we omit the dependency of $N_{i,j}$ on $k.$) 
Consequently, to identify the patterns in the coefficients of these identities as well as to provide closed formulas, it suffices to give explicit  
functional form to $N_{i,j},$ $D_k,$ and the constant term.
For $k=2,3,4,5,6$ the denominators on the right-hand side of the examples which occur  outside the big
parentheses are 2,4,6,8,10 leading to the conjecture that
\begin{equation}\label{dk}
	D_k = 2(k-1), \qquad k \ge 2.
\end{equation}

The obvious approach to finding the pattern in the $N_{i,j},$ creating a table whose rows
are labeled with $k$ and whose columns are labeled with $F^{(k)}_{m+i} F^{(k)}_{m+j}$ does not 
immediately yield results. The correct procedure that 
facilitates identification of patterns in the coefficients is to consider
the diagonal elements (that is, the case $i=j$) and the non-diagonal elements separately and in different ways.
The table of coefficients for the diagonal elements is presented in Table \ref{tnii}.

\begin{table}[H]
 \begin{small}
\begin{center}
\caption
{Coefficients, $N_{i,i}$  of $G_{m+i}^2, 0 \le i \le k-1$ for $2 \le k \le 6,$
based on \eqref{ex2}-\eqref{ex6}.}
\label{tnii}
{
\renewcommand{\arraystretch}{1.3}
\begin{tabular}{||c||c|c|c|c|c|c||} 
\hline \hline

$ $&$i=0$&$i=1$&$i=2$&$i=3$&$i=4$&$i=5$\\
\hline
\;&$G_m^2$&$G_{m+1}^2$&$G_{m+2}^2$&$G_{m+3}^2$&$G_{m+4}^2$&$G_{m+5}^2$\\
\hline
$k =2$&$0$&$0$&\;&\;&\;&\;\\
$k=3$&$-1$&$-4$&$-1$&\;&\;&\;\\
$k=4$&$-2$&$-8$&$-6$&$-2$&\;&\;\\
$k=5$&$-3$&$-12$&$-11$&$-8$&$-3$&\;\\
$k=6$&$-4$&$-16$&$-16$&$-14$&$-10$&$-4$\\

\hline \hline
\end{tabular}

}
\end{center}
 \end{small} 
\end{table}

This table shows clear linear patterns in each column
which naturally leads to a conjecture on the functional form
of the $N_{i,i}.$  The cases 
$i=0$ and $i \ge 1$ must be
treated separately.
For $i=0,$
\begin{equation}\label{n00}
	N_{0,0} = -(k-2),
\end{equation}
 while for $i \ge 1,$
\begin{equation}\label{nii}
	N_{i,i} = 4-(i+3)(k-i), \qquad k \ge 2, 1 \le  i \le k-1.
\end{equation}
Table \ref{tnii} has interesting symmetries as shown in \cite[A343125]{OEIS}.

The non-diagonal coefficients $N_{i,j}, i \neq j,$
 naturally form a 
3-dimensional solid rather than a triangle; more specifically for each $k,$ 
the coefficients $N_{i,j}, j \ge i+1,$ form a
triangle. For $k=6$ this triangle is shown in Table \ref{tiltj}.
For $2 \le k \le 5,$ the corresponding triangles are easy to construct;
(they may also be found (albeit crossed out) at 
\cite[Sequence A342955, history, item \#43]{OEIS}).

\begin{table}[H]
%\begin{tiny}
\begin{center}
\caption{Coefficients $N_{i,j}, 0 \le i \le k-2, i+1 \le j \le k-1,$ for the case $k=6,$
based on \eqref{ex6}. }
\label{tiltj}
{
\renewcommand{\arraystretch}{1.3}
\begin{tabular}{||c||c|c|c|c|c|c||} 
\hline \hline

$k=6$&$j=0$&$j=1$&$j=2$&$j=3$&$j=4$&$j=5$\\
\hline
$i=0$&\;&$-6$&$-4$&$-2$&$0$&$2$\\
$i=1$&\;&\;&$-8$&$-4$&$0$&$4$\\
$i=2$&\;&\;&\;&$-6$&$0$&$6$\\
$i=3$&\;&\;&\;&\;&$0$&$8$\\
$i=4$&\;&\;&\;&\;&\;&$10$\\ 
\hline \hline

\end{tabular}
}
\end{center}
%\end{tiny} 
\end{table}

Table \ref{tiltj}, like Table \ref{tnii}, exhibits linear patterns
which motivate the following conjecture on the functional form of
$N_{i,j}.$
\begin{equation}\label{niltj}
	N_{i,j}=2(i+1)(j-(k-2)), \qquad \text{for $0 \le i \le k-2, i+1 \le j \le k-1$}.
\end{equation} 

To clarify the flow of logic, we regard \eqref{dk}-\eqref{niltj} 
as definitions of $D_k, N_{i,i},$ and $N_{i,j}.$ The Main
Theorem will then prove that these definitions correctly provide
closed formulas for \eqref{TheProblem}.

\section{Main Theorem: Statement and Proof Outline}\label{sMT}

The previous section motivated the various components of the Main Theorem which we now state.

\begin{theorem}[Main Theorem]\label{maintheorem}
Fix $k \ge 2.$ Using \eqref{recursion},
let $\{G_n\}_{m \ge 0} =\{F^{(k)}_n\}_{n \ge 0}.$
Then, using  \eqref{dk}-\eqref{niltj},
\begin{equation}\label{main1}
	\displaystyle \sum_{i=0}^m G_i^2 = 
	\sum_{\substack {0 \le i \le k-1 \\ i \le j \le k-1}} 
	\frac{N_{i,j}}{D_k} G_{m+i} G_{m+j} - N_{k-1,k-1}.
\end{equation}
\end{theorem} 

\begin{comment} If the initial conditions 
in \eqref{recursion} are changed, then the theorem remains true
provided the $N_{k-1,k-1}$ summand is appropriately replaced.
A comment  after \eqref{main3}, will clarify this further.
\end{comment}

To prove Theorem \ref{maintheorem}, we first make some routine simplifications.

First, clearing denominators in \eqref{main1}, we obtain
\begin{equation}\label{main2}
	D_k \displaystyle \sum_{i=0}^m G_i^2 = 
\sum_{\substack {0 \le i \le k-1 \\ i \le j \le k-1}}  
 N_{i,j} G_{m+i} G_{m+j}-D_k N_{k-1,k-1}.
\end{equation}

If we replace $m$ by $m-1,$ in \eqref{main2}, we obtain
$$
	D_k \displaystyle \sum_{i=0}^{m-1} G_i^2 = 
	\sum_{\substack {0 \le i \le k-1 \\ i \le j \le k-1}} 
	N_{i,j} G_{m-1+i} G_{m-1+j}-D_k N_{k-1,k-1}. 
$$
If we now take the difference of these last two equations,
 we see that to prove \eqref{main1}, it suffices to prove 
 \begin{align}\label{main3}
	D_k G_m^2 &= 
	\sum_{\substack {0 \le i \le k-1 \\ i \le j \le k-1}}  N_{i,j}
	\biggl( G_{m+i} G_{m+j} &- G_{m-1+i} G_{m-1+j} \biggr) \\
	&=\sum_{\substack {0 \le i \le k-1 \\ i \le j \le k-1}}  
	{N}_{i,j}G_{m+i} G_{m+j}
	 &- \sum_{\substack {-1 \le i \le k-2 \\ i \le j \le k-2}}  
	{N}_{i+1,j+1} G_{m+i} G_{m+j}   \nonumber  \\
	 &= \sum_{\substack {-1 \le i \le k-2 \\ i \le j \le k-2}}  
	{N'}_{i,j} G_{m+i} G_{m+j}, 					\nonumber 
\end{align} 
where  i) summands with a $G_{m+k-1}$ multiplicand have been eliminated using
\eqref{recursion} which implies that $G_{m+k-1}= \sum_{i=-1}^{k-2} G_{m+i},$
and ii)   ${N'}_{i,j}$ is the resulting linear combination of the $N_{i,j}.$

\begin{comment} Notice how the $N_{k-1,k-1} D_k$ term vanished because of  the subtraction.
This provides a basis for proving the Main Theorem when different initial conditions
are used. \end{comment} 

In the sequel, we will abuse language and refer to the ``sides" of \eqref{main3}
referring to the items that are set
equal in \eqref{main3}. We similarly will refer to 
i) $\sum_{\substack {0 \le i \le k-1 \\ i \le j \le k-1}}      {N}_{i,j}G_{m+i} G_{m+j}$
and ii)
	 $- \sum_{\substack {-1 \le i \le k-2 \\ i \le j \le k-2}}  {N}_{i+1,j+1} G_{m+i} G_{m+j}$, 
as \emph{the two summands on the middle 
side of \eqref{main3}}. This terminology, while slightly non-standard, should cause no confusion.
  
We next outline the key idea in the proof of \eqref{main3}.

By \eqref{n00}-\eqref{niltj}, $N_{i,j}$ is a  polynomial in at most three
variables $i,j$ and $k,$ of degree at most 2, and consequently, ${N'}_{i,j}$ which is some linear
combination of the $N_{i,j}$ is also a polynomial in at most 3 variables
of degree at most 2. 

Therefore, to prove \eqref{main1}, for which it suffices to prove \eqref{main3}, 
it suffices to i) calculate the ${N'}_{i,j}, -1 \le i \le k-2, i \le j \le k-2$
and then ii) \emph{algebraically verify}  that the coefficients of $G_{m+i} G_{m+j}$ 
on each side of 
\eqref{main3} are identical, or equivalently, that 
iii) ${N'}_{i,j}=D_k$ if $i=0=j$ 
and 0 otherwise.

\section{The Seven Coefficient Groups} \label{s7}

Before presenting the details of the proof, we have one
subtlety to deal with. The algebraic verification method
requires checking polynomial-coefficient equality over all
index pairs $(i,j).$ But these index pairs lie in the upper
half of the  $k \times k$ square; hence, the number of verifications
could be going to infinity.

It turns out that for any $k \ge 2,$ there are only seven distinct sets of index pairs that 
need to be considered for the proof of the Main Theorem.
The polynomial
form of ${N'}_{i,j}$ is identical for each of these seven sets. Hence,
we only need algebraically verify at most seven polynomial equalities.
 
We will let the letters, $A, \dotsc, G$ 
indicate these seven sets. The seven sets of index pairs
are formally defined as follows.
 
\begin{align*} 
	A  &= \{i=-1, j=-1\} \\
	B  &= \{i=0, j=0\} \\
	C  &= \{1 \le i \le k-2, j=i\} \\
	D  &= \{i=-1, 0 \le j \le k-3\} \\
	E  &= \{i=-1, j=k-2\} \\
	F  &= \{0 \le i \le k-4, i+1 \le j \le k-3\} \\
	G  &= \{0 \le i \le k-3, j=k-2\}.  
\end{align*} 
It is straightforward to check that these 7 sets are mutually 
exclusive and completely cover the set of index pairs
 $\{(i,j): -1 \le i \le k-2, i \le j \le k-2\}.$ Table \ref{tABC}
corresponds to the formal definition just given.

In the sequel, instead of using this formal definition
we will use \emph{label indexed squares} 
which visually depict these seven groups.  Labeled index squares
are a contribution of this paper; they facilitate following the flow of
a complex proof.

\begin{table}[H]
%\begin{tiny}
\begin{center}
\caption{The seven groups of $G_{m+i} G_{m+j}$}
\label{tABC}
{
\renewcommand{\arraystretch}{1.3}
\begin{tabular}{||c||c|c|c|c|c|c|c|c||} 
\hline \hline
\;&$G_{m-1}$&$G_m$&$G_{m+1}$&$\dotsi$&$G_{m-(k-4)}$&$G_{m-(k-3)} $&$G_{m-(k-2)}$\\
\hline \hline
$G_{m-1}$&$A$&$D$&$\dotsi$&$\dotsi$&$\dotsi$&$D$&$E$\\
$G_m$&\;&$B$&$F$&$\dotsi$&$\dotsi$&$F$&$G$\\
$G_{m+1}$&\;&\;&$C$&$\ddots$&$\ddots$&$\vdots$&$\vdots$\\
$\vdots$&\;&\;&\;&$\ddots$&$\ddots$&$\vdots$&$\vdots$\\
$G_{m-(k-4)}$&\;&\;&\;&\;&$\ddots$&$F$&$\vdots$\\
$G_{m-(k-3}) $&\;&\;&\;&\;&\;&$\ddots$&$G$\\
$G_{m-(k-2)}$&\;&\;&\;&\;&\;&\;&$C$\\
\hline \hline

\end{tabular}
}
\end{center}
%\end{tiny}
 \end{table}

The following examples are illustrative of the use of the index-pair sets
and labeled index squares.

\begin{example} $ \sum_{A} N_{i+1,j+1} G_{m+i} G_{m+j} = N_{0,0} G_{m-1}^2.$
In this example, by  Table \ref{tABC},
$A= \{i=-1,j=-1\}.$ Hence $N_{i+1,j+1} G_{m+i} G_{m+j} = N_{0,0} G_{m-1}^2$ as 
required. \end{example}

\begin{example} A set of letters separated by commas will indicate unions
of sets, \\
 for example, $\sum_{A,B,C} N_{i,j} G_{m+i} G_{m+j} 
= \sum_{-1 \le i \le k-2} N_{i,i} G_{m+i}^2.$
 \end{example}

\begin{example} 
$\sum_{\substack{0 \le i \le k-2 \\ 0 \le j \le k-2 \\ i \neq j}} N_{i,k-1} G_{m+i} G_{m+j}= \sum_{F,G} \bigg( N_{i,k-1} + N_{j,k-1}\bigg)G_{m+i} G_{m+j}.$  
To clarify the subtleties in this sum,  consider, by way of illustration
two index pairs, say $i=2,j=4$ and $i=4,j=2.$ 
For $i=2,j=4$ the corresponding summand is 
$N_{2,k-1} G_{m+2} G_{m+4}= N_{i,k-1} G_{m+i}G_{m+j}.$
However, since \eqref{main3} is summed over pairs of indices with $i \le j,$ it follows that
for $i=4,j=2$ the corresponding summand is
$N_{4,k-1} G_{m+2} G_{m+4}=  N_{j,k-1} G_{m+i}G_{m+j}.$ 
As this example shows, the sum goes over both the index pairs with $i<j$ as well as over the index pairs with $j<i.$
\end{example}

\section{Detailed Outline of Proof}\label{sP}

A general overview of the proof begun
in Section \ref{sMT}, shows that the proof of the Main Theorem 
is accomplished through
verification. The actual verification will involve distributing 24 polynomial-summands  
over the seven index-sets and then verifying the polynomial-equality corresponding to
each index-set on 
each side of \eqref{main3}. To provide complete details
of the flow
of logic, this section organizes the proof as a series of enumerated steps with accompanying appropriate
references to equations and tables.

\begin{enumerate}[Step i.]
\item We start with \eqref{main1} the theorem statement to be proven, the Main Theorem.
\item We then clear denominators, 
reducing proof of \eqref{main1} to \eqref{main2}
\item Equation \eqref{main3}   
accomplishes three things. First, it eliminates the summation
of the $G_i^2$ (and replaces it with a single $G_m^2.$) Additionally, this step replaces $G_{m+k-1}$ with
$\sum_{i=-1}^{k-2} G_{m+i}$ using \eqref{recursion}.
\item Equation \eqref{main3} also introduces the ${N'}_{i,j},$ 
linear combinations of the $N_{i,j},$
and shows that to prove the Main Theorem it suffices to 
prove the equality of polynomial coefficients on all sides 
of \eqref{main3}.  
\item Equations \eqref{summand123} and \eqref{summand45} break down 
the two summands on the middle
side of \eqref{main3}
into three and two summands respectively called Summand-I through Summand-III
and Summand-IV through Summand-V. 
\item Tables \ref{tSummandI} through \ref{tSummandV}  
analyze these fives sums in terms of the index pairs over which they
are defined and the corresponding $N_{i,j}$ coefficients.
This analysis is facilitated by use of   labeled index squares.
\item It immediately follows, that if for each index set, $A,\dotsc,G,$ we sum over the five summands listed in Tables \ref{tSummandI} 
through \ref{tSummandV} we will obtain the polynomial form of
${N'}_{i,j}$ of \eqref{main3}. 
\item Tables \ref{tSummary} and \ref{tSummarynumber}
accomplish the summation mentioned in Step vii, by
first presenting a $7 \times 5$ table listing the coefficient forms of the $N_{i,j}$
over each of the seven index sets (the rows) 
and each of the five summands (the columns).
\item Then, using \eqref{dk}-\eqref{niltj}, Table \ref{tSummarynumber} evaluates each $N_{i,j}$ entry in
Table \ref{tSummary} and additionally sums the polynomials across rows. These sums are
the ${N'}_{i,j}$ of \eqref{main3}. They are linear combinations of the $N_{i,j}.$ 
The summations are accomplished through
pencil-and-paper manipulation.
\item A glance at this sum column in Table \ref{tSummarynumber} 
confirms equality  of the sides in \eqref{main3}. More specifically,
  for $i=0=j,$  the polynomial coefficient is $D_k$ and 0 otherwise.
\item This completes the proof
\end{enumerate} 

\textbf{Illustration of the proof on one summand.}
To further clarify these 11 steps we illustrate them with
the proof of the equality of polynomial coefficients of $G_{m-1}G_{m-1}$ on all sides of \eqref{main3}. 

\begin{proof} 
There are two summands on the middle side of 
\eqref{main3}. 

\emph{The first summand.} When 
$i=k-1=j$ we obtain the summand 
$N_{k-1,k-1} G_{m+k-1}^2.$ If we then expand
$G_{m+k-1}^2$ by applying \eqref{recursion} we obtain
$N_{k-1,k-1} \biggl( \sum_{i=-1}^{k-2} G_{m+i}\biggr)^2.$
Thus we have a contribution of $N_{k-1,k-1} G_{m-1}^2.$

\emph{The second summand.} When $i=-1=j$
we obtain the summand $N_{0,0} G_{m-1}^2.$
 
Thus the total contribution of the middle side of \eqref{main3} to summands involving $G_{m-1}^2$ is $(N_{-,1-1}+N_{k-1,k-1}) G_{m-1}^2.$

This implies that ${N'}_{-1,-1} = N_{k-1,k-1} - N_{0,0}.$ 
By \eqref{n00}-\eqref{nii},  ${N'}_{-1,-1}=0$ as required.
\end{proof}

This completes the proof of one of the seven required
verifications of equality of polynomial coefficients
on all sides of \eqref{main3}.

 A similar argument applies to the other six required verifications.

We can use this sample proof to illustrate the steps enumerated above.
 Step iv states the overall goal of proving the equality of coefficients of $G_{m+i} G_{m+j}$
on all sides of \eqref{main3}.
Step v, identifying the contribution of $N_{k-1,k-1}$
to ${N'}_{-1,-1}$ is found in Summand-III of \eqref{summand123}. Step vi further analyzes Summand-III
with a labeled index-square in Table \ref{tSummandIII}
where we find a coefficient of $N_{k-1,k-1}$ associated with index set $A$
as required. Step viii, summarizes this association of $N_{k-1,k-1}$ 
with index set $A$ and Summand-III  
in Table \ref{tSummary} in the row labeled A and the 
column labeled Summand-III. Then Step ix, 
evaluates this $N_{k-1,k-1}$ cell entry as a polynomial
using \eqref{n00}-\eqref{nii}. This evaluation may be found
in Table \ref{tSummarynumber} in the row labeled A and the 
column labeled Summand-III.

A similar argument, or trace of logical flow, applies to the
contribution of $N_{0,0}$ to ${N'}_{-1,-1}.$

Finally, in Step x, 
Table \ref{tSummarynumber} evaluates and sums the 
polynomials  in its first
row showing that ${N'}_{-1,-1} = 0.$ This
provides one of the seven required
verifications needed for the entire proof which when done
accomplishes completion of the proof as indicated in Step xi.

In the sequel, we will refer to the steps in this section as we implement each one. Steps i - iv have already been  accomplished.

\section{The Five Summands}\label{s5s}

This section implements Step v of the proof, breaking down the two summands on the middle side of
\eqref{main3} into five summands called Summand-I through Summand-V.
Corresponding to the summands on the middle side of \eqref{main3} we have the
following decompositions.

For the first summand on the middle side of \eqref{main3}, we have
\begin{multline}\label{summand123}
\sum_{\substack {0 \le i \le k-1 \\ i \le j \le k-1}}  N_{i,j}G_{m+i} G_{m+j}
 	= \text{Summand-I}+ \text{Summand-II} + \text{Summand-III} \\
	=\sum_{\substack {0 \le i \le k-2 \\ i \le j \le {k-2}}}  N_{i,j}G_{m+i} G_{m+j}+
	\sum_{\substack {0 \le i \le k-2 \\ j = k-1}}  N_{i,j}G_{m+i} 
	\biggl(\sum_{i=-1}^{k-2} G_{m+i}\biggr)  +
	  N_{k-1,k-1}\biggl(\sum_{i=-1}^{k-2} G_{m+i}\biggr)^2,
\end{multline}
where in summand-II and summand-III we replace 
$G_{m+k-1}$ with $\sum_{i=-1}^{k-2} G_{m+i}$ using \eqref{recursion}.

For the second summand on the middle side of \eqref{main3}, we have
\begin{multline}\label{summand45}
- \sum_{\substack {-1 \le i \le k-2 \\ i \le j \le k-2}}  
	N_{i+1,j+1} G_{m+i} G_{m+j} =\text{Summand-IV} + \text{Summand-V} \\
	=-\sum_{-1 \le i \le k-2}  N_{i+1,i+1}G_{m+i}^2  -
	\sum_{\substack {-1 \le i \le k-3 \\ i <j \le k-2}}  N_{i+1,j+1}G_{m+i} G_{m+j}.
\end{multline}

\section{Index tables for Summands-I through Summands-V}\label{s5}  

This section implements Step vi of the proof, analyzing the five summands of the last 
section according to the index pairs over which they are summed. For the convenience
of the reader, the five summands are repeated in the header of each table. The labeled
index squares provide visual aids and review the definition of the seven index sets.
Each cell of these labeled index squares contain i) the index letter for that pair of indices
and ii) the linear combination of $N_{i,j}$ for that cell. These two items are separated by a comma.
The index-squares should be fairly straightforward to read and check.

\begin{table}[H]
%\begin{tiny}
\begin{center}
\caption{$\text{Summand-I}=  
\sum_{\protect\substack {0 \le i \le k-2 \\ i \le j \le k-2}}  N_{i,j} G_{m+i} G_{m+j}
=\sum_{B,C,F,G} N_{i,j} G_{m+i} G_{m+j}$} 
\label{tSummandI}
{
\renewcommand{\arraystretch}{1.3}
\begin{tabular}{||c||c|c|c|c|c|c|c|c||} 

 \hline \hline

\;&$G_{m-1}$&$G_m$&$G_{m+1}$&$\dotsi$&$G_{m-(k-4)}$&$G_{m-(k-3}) $&$G_{m-(k-2)}$\\ 
\hline $G_{m-1}$&$ $&$ $&$ $&$ $&$ $&$ $&$ $\\
\hline $G_m$&\;&$B, N_{0,0}$&$F, N_{i,j}$&$\dotsi$&$\dotsi$&$F, N_{i,j}$&$G, N_{i,k-2}$\\ 
\hline $G_{m+1}$&\;&\;&$C, N_{i,i}$&$\ddots$&$\ddots$&$\vdots$&$\vdots$\\
\hline $\vdots$&\;&\;&\;&$\ddots$&$\ddots$&$\vdots$&$\vdots$\\
\hline $G_{m-(k-4)}$&\;&\;&\;&\;&$\ddots$&$F, N_{i,j}$&$\vdots$\\
\hline $G_{m-(k-3}) $&\;&\;&\;&\;&\;&$\ddots$&$G, N_{i,k-2}$\\
\hline $G_{m-(k-2)}$&\;&\;&\;&\;&\;&\;&$C, N_{i,i}$\\

\hline \hline

\end{tabular}
}
\end{center}

%\end{tiny}
 \end{table}
  
\begin{table}[H]
%\begin{tiny}
\begin{center}

\caption
{
$\text{Summand-II}=  
\sum_{\protect\substack {0 \le i \le {k-2} \\ j = k-1}}  N_{i,j}G_{m+i} 
	\biggl(\sum_{i=-1}^{k-2} G_{m+i}\biggr) 
=
\sum_{B,C,F,G} N_{i,k-1} G_{m+i} G_{m+j}+
\sum_{D,E,F,G} N_{j,k-1} G_{m+i} G_{m+j}$
} 
\label{tSummandII}
{
\renewcommand{\arraystretch}{1.3}
\begin{tabular}{||c||c|c|c|c|c|c|c|c||} 

\hline \hline

\;&$G_{m-1}$&$G_m$&$G_{m+1}$&$\dotsi$&$G_{m-(k-4)}$&$G_{m-(k-3}) $&$G_{m-(k-2)}$\\
\hline $G_{m-1}$&$ $&$D, N_{j,k-1}$&$\dotsc$&$\dotsc$&$\dotsc$&$D, N_{j,k-1}$&$E, N_{k-2, k-1}$\\
\hline $G_m$&\;&$B, N_{0,k-1}$&$F, N_{i,k-1} + N_{j, k-1}$&$\dotsi$&$\dotsi$&$F, N_{i,k-1} + N_{j, k-1}$&$G, N_{i,k-1} + N_{k-2,k-1}$\\
\hline $G_{m+1}$&\;&\;&$C, N_{i,k-1}$&$\ddots$&$\ddots$&$\vdots$&$\vdots$\\
\hline $\vdots$&\;&\;&\;&$\ddots$&$\ddots$&$\vdots$&$\vdots$\\
\hline $G_{m-(k-4)}$&\;&\;&\;&\;&$\ddots$&$F, N_{i,k-1} + N_{j, k-1}$&$\vdots$\\
\hline $G_{m-(k-3}) $&\;&\;&\;&\;&\;&$\ddots$&$G, N_{i,k-1} + N_{k-2,k-1}$\\
\hline $G_{m-(k-2)}$&\;&\;&\;&\;&\;&\;&$C, N_{i,k-1}$\\

\hline \hline

\end{tabular}
}
\end{center}

%\end{tiny}
 \end{table}

\begin{table}[H]
%\begin{tiny}
 \begin{center}

\caption
{
$\text{Summand-III}=  
	  N_{k-1,k-1}\biggl(\sum_{i=-1}^{k-2} G_{m+i} \biggr)^2
=
\sum_{A,B,C} N_{k-1,k-1} G_{m+i}^2 +
\sum_{D,E,F,G} 2N_{k-1,k-1} G_{m+i}G_{m+j} 
$
} 
\label{tSummandIII}
{
\renewcommand{\arraystretch}{1.3}
\begin{tabular}{||c||c|c|c|c|c|c|c|c||} 

\hline \hline

\;&$G_{m-1}$&$G_m$&$G_{m+1}$&$\dotsi$&$G_{m-(k-4)}$&$G_{m-(k-3}) $&$G_{m-(k-2)}$\\
\hline $G_{m-1}$&$A, N_{k-1,k-1}$&$D,  2N_{k-1,k-1}$&$\dotsc$&$\dotsc$&$\dotsc$&$D,  2N_{k-1,k-1}$&$E,  2N_{k-1,k-1}$\\
\hline $G_m$&\;&$B, N_{k-1,k-1}$&$F, 2N_{k-1,k-1}$&$\dotsi$&$\dotsi$&$F, 2N_{k-1,k-1}$&$G,  2N_{k-1,k-1}$\\
\hline $G_{m+1}$&\;&\;&$C, N_{k-1,k-1}$&$\ddots$&$\ddots$&$\vdots$&$\vdots$\\
\hline $\vdots$&\;&\;&\;&$\ddots$&$\ddots$&$\vdots$&$\vdots$\\
\hline $G_{m-(k-4)}$&\;&\;&\;&\;&$\ddots$&$F, 2N_{k-1,k-1}$&$\vdots$\\
\hline $G_{m-(k-3}) $&\;&\;&\;&\;&\;&$\ddots$&$G,  2N_{k-1,k-1}$\\
\hline $G_{m-(k-2)}$&\;&\;&\;&\;&\;&\;&$C, N_{k-1,k-1} $\\

\hline \hline

\end{tabular}
}
\end{center}

%\end{tiny}
 \end{table}

\begin{table}[H]
%\begin{tiny}
\begin{center}

\caption
{
$\text{Summand-IV}=  
 -\sum_{-1 \le i \le k-2}  N_{i+1,i+1}G_{m+i}^2 
=
-\sum_{A, B, C} N_{i+1,i+1} G_{m+i}^2$
} 
\label{tSummandIV}
{
\renewcommand{\arraystretch}{1.3}
\begin{tabular}{||c||c|c|c|c|c|c|c|c||} 

\hline \hline

\;&$G_{m-1}$&$G_m$&$G_{m+1}$&$\dotsi$&$G_{m-(k-4)}$&$G_{m-(k-3}) $&$G_{m-(k-2)}$\\
\hline $G_{m-1}$&$A, -N_{0,0}$&$ $&$ $&$ $&$ $&$ $&$ $\\
\hline $G_m$&\;&$B, -N_{1,1}$&$ $&$ $&$ $&$ $&$ $\\
\hline $G_{m+1}$&\;&\;&$C,  -N_{i+1,i+1}$&$ $&$ $&$ $&$ $\\
\hline $\vdots$&\;&\;&\;&$\ddots$&$ $&$ $&$ $\\
\hline $G_{m-(k-4)}$&\;&\;&\;&\;&$\ddots$&$ $&$ $\\
\hline $G_{m-(k-3}) $&\;&\;&\;&\;&\;&$\ddots$&$ $\\
\hline $G_{m-(k-2)}$&\;&\;&\;&\;&\;&\;&$C,  -N_{i+1,i+1}$\\

\hline \hline

\end{tabular}
}
\end{center}

%\end{tiny}
 \end{table}

\begin{table}[H]
%\begin{tiny}
\begin{center}

\caption
{
$\text{Summand-V}=  
-\sum_{\protect\substack {-1 \le i \le k-3 \\ i <j \le k-2}}  N_{i+1,j+1}G_{m+i} G_{m+j}
=
-\sum_{D, E ,F, G} N_{i+1,j+1} G_{m+i} G_{m+j}$
} 
\label{tSummandV}
{
\renewcommand{\arraystretch}{1.3}
\begin{tabular}{||c||c|c|c|c|c|c|c|c||} 

\hline \hline

\;&$G_{m-1}$&$G_m$&$G_{m+1}$&$\dotsi$&$G_{m-(k-4)}$&$G_{m-(k-3}) $&$G_{m-(k-2)}$\\
\hline $G_{m-1}$&$ $&$D,  -N_{0,j+1}$&$\dotsc$&$\dotsc$&$\dotsc$&$D,  -N_{0,j+1}$&$E, -N_{0,k-1}$\\
\hline $G_m$&\;&$ $&$F, -N_{i+1, j+1}$&$\dotsc$&$\dotsc$&$F,  -N_{i+1, j+1}$&$G, -N_{i+1,k-1}$\\
\hline $G_{m+1}$&\;&\;&$ $&$\ddots$&$\ddots$&$\vdots$&$\vdots$\\
\hline $\vdots$&\;&\;&\;&$ $&$\ddots$&$\vdots$&$\vdots$\\
\hline $G_{m-(k-4)}$&\;&\;&\;&$ $&$ $&$F,  -N_{i+1, j+1}$&$\vdots$\\
\hline $G_{m-(k-3}) $&\;&\;&\;&\;&\;&$ $&$G, -N_{i+1,k-1}$\\
\hline $G_{m-(k-2)}$&\;&\;&\;&\;&\;&\;&$ $\\

\hline \hline

\end{tabular}
}
\end{center}

%\end{tiny}
 \end{table}

\section{Distribution of 
\texorpdfstring{$N_{i,j}$}{N(i,j)} Across the Index-Sets}\label{sDistribution}

This section implements Steps vii - xi of the proof outline. Step vii states the goal,
to show that the ${N'}_{i,j}$ are polynomially equal on all  sides of \eqref{main3}.
  
To accomplish this goal, first, in Step viii, Table \ref{tSummary}
summarizes the coefficients for each of the five summands of the previous section for each
index set. The coefficients are expressed in terms of $N_{i,j}.$

\begin{table}[H]
\begin{center}

\caption{Coefficient of $G_{m+i} G_{m+j}$ in  Summand-I through Summand-V as given 
in \eqref{summand123}-\eqref{summand45}.} 
\label{tSummary}
 
{
\renewcommand{\arraystretch}{1.3}

\begin{tabular}{||c|c|c|c|c|c||} 

\hline \hline
%\begin{tiny}
 
\;&$Summand-I$&$Summand-II$&$Summand-III$&$Summand-IV$&$Summand-V$\\
\hline $A$&\;&\;&$N_{k-1, k-1}$&$-N_{0,0}$&\;\\
\hline $B$&$N_{0,0}$&$N_{0,k-1}$&$N_{k-1, k-1}$&$-N_{1,1}$&\;\\
\hline $C$&$N_{i,i}$&$N_{i, k-1}$&$N_{k-1, k-1}$&$-N_{i+1,i+1}$&\;\\
\hline $D$&\;&$N_{j, k-1}$&$2N_{k-1,k-1}$&\;&$-N_{0,j+1}$\\
\hline $E$&\;&$N_{k-2,k-1}$&$2N_{k-1,k-1}$&\;&$-N_{0, k-1}$\\
\hline $F$&$N_{i, j}$&$N_{i,k-1}+N_{j,k-1}$&$2N_{k-1,k-1}$&\;&$-N_{i+1, j+1}$\\
\hline $G$&$N_{i, k-2}$&$N_{i,k-1}+N_{k-2,k-1}$&$2N_{k-1,k-1}$&\;&$-N_{i+1,k-1}$\\

\hline \hline
 %\end{tiny}
\end{tabular}
}

\end{center}

 \end{table}

Next, we implement Step ix. Using \eqref{dk}-\eqref{niltj} we evaluate each entry in Table \ref{tSummary}
as a polynomial in the variables $i,j,k.$

\begin{table}[H]
\begin{center}

\caption{Coefficient of $G_{m+i} G_{m+j}$ in equations Summand-I through Summand-V.} 
\label{tSummarynumber}
% \begin{tiny}
{
\renewcommand{\arraystretch}{1.3}
\resizebox{\textwidth}{!}{
\begin{tabular}{||c|c|c|c|c|c|c||} 
\hline \hline
 
\;&$Summand-I$&$Summand-II$&$Summand-III$&$Summand-IV$&$Summand-V$&$Sum ={N'}_{i,j}$\\
\hline $A$&\;&\;&$2 - k$&$k - 2$&\;&$0$\\
\hline $B$&$2 - k$&$2$&$2 - k$&$4(k-1) - 4$&\;&$2k-2=D_k$\\
\hline $C$&$4 - (i + 3)(k - i)$&$2(i + 1)$&$2 - k$&$(i + 4)(k - i - 1) - 4$&\;&$0$\\
    \hline $D$&\;&$2(j + 1)$&$4 - 2k$&\;&$-2(j - (k-3))$&$0$\\
\hline $E$&\;&$2(k - 1)$&$4 - 2k$&\;&$-2$&$0$\\
\hline $F$&$2(i+1)(j-(k-2))$&$2(i + 1)  + 2(j + 1)$&$4 - 2k$&\;&$-2(i + 2)(j - (k - 3))$&$0$\\
\hline $G$&$0$&$2(i + 1) + 2(k - 1)$&$4 - 2k$&\;&$-2(i + 2)$&$0$\\

\hline \hline
\end{tabular}}
}
\end{center}

 %\end{tiny} 
 \end{table}

The sum column, accomplished through paper and pencil summation of polynomials of each row, implements
Step x of the proof outline and calculates the $N'_{i,j}.$ As can be seen from \eqref{main3},
the sum column confirms the algebraic verification of polynomial equality 
of the sides of  \eqref{main3} since, by \eqref{dk}, ${N'}_{i,j}$ is equal to $D_k$ for $i=0=j$ and 0 otherwise as required. 

As observed in Step xi, the proof of the Main Theorem is complete.

\section{Conclusion}\label{sC}

This paper has proven the Main Theorem identifying the patterns in the coefficients in the
closed formulas for the sums of squares of $k$-bonacci numbers. Along the way, the paper
presented the algebraic verification method and introduced the labeled index-squares as a visual
aid facilitating the proof. 

This result and other recent results suggest a possible new trend in approaching identities. \cite{HendelR}
reviews several stages (not necessarily consecutive) in the history of (Fibonacci-Lucas) identities. One stage is simply a concern for punchy, cute, unexpected equalities. Another stage is a concern
for methods, such as Binet form, matrices, generating functions etc. 
The goal of this stage 
is to see how existing identities can be proven with a given method as well as whether 
new methods provide more elegant proofs. Still another stage is generalizing known Fibonacci-Lucas
identities, for example, Cassini, Catalan, d'Ocagne, to other second order recursions.

Recently however, a new stage has emerged, families of identities.
For example,  Melham \cite{Melham} proved a remarkable 
generalization of an identity of Aurifeuille for a family of identities where
the number of summands is going to infinity. He used the Dresel Verification Theorem
\cite{Dresel} to prove several individual cases thus creating credibility to the result
which is still open. Hendel has proven results about families of identities where the
order of the recursions are going to infinity. In \cite{HendelTT}, the family of identities
corresponds to a Taylor series with the Taylor polynomials of this Taylor series equaling
the characteristic polynomials which define the recursions in this family. In \cite{HendelS},
a result about a family of recursions whose orders are going to infinity, is proven uniformly for all orders,
by examining the divisibility properties of characteristic polynomials of the minimal recursions
of the family members. In this paper, a result about sums of squares, with initial cases each
requiring different methods and ingenuity, was proven for a family of recursions, whose orders
are going to infinity by a simple (but tedious) algebraic verifciation of polynomial-equality.

We believe these results and methods could point to a new trend 
that emphasizes studying families of identities.
We encourage other researchers to use the
methods mentioned to achieve new results.

%\begin{acknowledgements}
%\end{acknowledgements}



\begin{thebibliography}{99}
 
 

\bibitem{Dresel}
L. A. G. Dresel,
``Transformations of Fibonacci Lucas Identities,"
\emph{Applications of Fibonacci Numbers}, 
\textbf{Vol. 5}, 169--184, 
Kluwer Acad. Publ., Dordrecht, 1993


\bibitem{HendelR}
Russell Jay Hendel;
``Review of Subsequences of Fibonacci and Lucas polynomials with geometric subscripts,  
Fibonacci Quart., 50,  27--35,"
\emph{Mathematical Reviews}, \textbf{\#2892007}, (2012).
 
 
\bibitem{HendelTT}
Russell Jay Hendel, 
``Recursive Triangles Appearing Embedded in Recursive Families,"
\emph{Fibonacci Quarterly},
\textbf{58}, (2020),
135--144.

\bibitem{HendelS}
Russell Jay Hendel,
`` A Method for Uniformly Proving a Family of Identities,"
 \emph{Fibonacci Quarterly}, To appear,2022.

 \bibitem{Howard_Cooper}
Fred T. Howard and Curtis Cooper, 
``Some Identities for $r$-Fibonacci Numbers,"
\emph{The Fibonacci Quarterly}, \textbf{49(3)}, 231--242, (2011).

\bibitem{Koshy}
T. Koshy,
\emph{Fibonacci and Lucas Numbers with Applications,}
John Wiley and Sons, New York, 2001.

\bibitem{Jakubczyk}
Zbigniew Jakubczyk, 
``Sums of Squares of Tribonacci Numbers, 
(Solution of Advanced Problem H-715),"
\emph{Fibonacci Quarterly}, \textbf{51(3)}, 285--286, (2013); "Alternating Sums of High Powers of Fibonacci Numbers, (Solution of Advanced Problem H-719)," \emph{Fibonacci Quarterly}, 
\textbf{51(4)}, 379--380, (2013).

 
\bibitem{Melham}
R. S. Melham,
``On a Classical Fibonacci Identity of
Aurifeuille",
\emph{Fibonacci Quarterly}, \textbf{54(1)},(2016), 19--22.

\bibitem{OEIS} \emph{The on-line Encyclopedia of Integer Sequences,}
http://oeis.org

\bibitem{Ohtsuka}
Hideyuki Ohtsuka,
``Advanced Problem H-715,"
\emph{Fibonacci Quarterly}, 50(1), 90, (2012).
 

 \bibitem{Prodinger_Selkirk}
Helmut Prodinger and Sarah J. Selkirk,
``Sums of Squares of Tetranacci Numbers: A Generating Function Approach,"
\emph{arXiv}, https://arxiv.org/abs/1906.08336, (2019). 


\bibitem{Schumacher3}
Raphael Schumacher,
``Explicit Formulas for Sums Involving the Squares
 of the First $n$ Tribonacci Numbers",
\emph{Fibonacci Quarterly}, \textbf{58},
194--202, (2020).

 
\bibitem{Schumacher4}
Raphael Schumacher,
``How to Sum the Squares of the Tetranacci Numbers
and the Fibonacci $m$-STEP Numbers,"
\emph{Fibonacci Quarterly}, 57(2), 168--175, (2019). 

 
 


% etc
\end{thebibliography}


\end{document}